\documentclass[12pt]{article}
\usepackage{mathrsfs}
\usepackage{amsfonts}
\usepackage{graphicx}
\usepackage{amsfonts, amsmath, amssymb, latexsym, epsfig}
\newcommand{\blst}{\begin{trivlist}}
\newcommand{\elst}{\end{trivlist}}

\newtheorem{thm}{Theorem}[section]
\newtheorem{prop}{Proposition}[section]
\newtheorem{cor}{Corollary}[section]
\newtheorem{lem}{Lemma}[section]
\newtheorem{conj}{Conjecture}[section]
\newtheorem{exa}{Example}[section]
\newtheorem{defn}{Definition}[section]

\newcommand{\ben}{\begin{enumerate}}
\newcommand{\een}{\end{enumerate}}
\newcommand{\ble}{\begin{lem}}
\newcommand{\ele}{\end{lem}}
\newcommand{\bth}{\begin{thm}}
\renewcommand{\eth}{\end{thm}}
\newcommand{\bpr}{\begin{prop}}
\newcommand{\epr}{\end{prop}}
\newcommand{\bco}{\begin{cor}}
\newcommand{\eco}{\end{cor}}
\newcommand{\bcon}{\begin{conj}}
\newcommand{\econ}{\end{conj}}
\newcommand{\bde}{\begin{defn}}
\newcommand{\ede}{\end{defn}}
\newcommand{\bex}{\begin{exa}}
\newcommand{\eex}{\end{exa}}
\newcommand{\barr}{\begin{array}}
\newcommand{\earr}{\end{array}}
\newcommand{\btab}{\begin{tabular}}
\newcommand{\etab}{\end{tabular}}
\newcommand{\beq}{\begin{equation}}
\newcommand{\eeq}{\end{equation}}

\newcommand{\bea}{\begin{eqnarray*}}
\newcommand{\eea}{\end{eqnarray*}}
\newcommand{\beaa}{\begin{eqnarray}}
\newcommand{\eeaa}{\end{eqnarray}}
\newcommand{\bce}{\begin{center}}
\newcommand{\ece}{\end{center}}
\newcommand{\bpi}{\begin{picture}}
\newcommand{\epi}{\end{picture}}
\newcommand{\bfi}{\begin{figure} \begin{center}}
\newcommand{\efi}{\end{center} \end{figure}}

\newcommand{\bsl}{\begin{slide}{}}
\newcommand{\esl}{\end{slide}}

{}
{}
{}
{}
{}
{}
{}
{}
{}
{}
{}
{}
{}
{}
{}
{}
{}
{}
{}
{}
{}
{}
{}
{}
{}
{}
{}
{}
\newenvironment{proof}{
\par
\noindent {\bf Proof.}\rm}{\mbox{}\hfill\rule{0.5em}{0.809em}\par}
\begin{document}
\title{Enumerations for Permutations by Circular Peak Sets}
\author{Pierre Bouchard$^{a}$
 \and Hungyung Chang$^{b}$
 \and Jun Ma$^{c,}$\thanks{Email address of the corresponding author: majun@math.sinica.edu.tw}\\
 \and Jean Yeh$^{d,}$\thanks{jean.yh@ms45.url.com.tw}}
\maketitle \vspace*{-1.2cm}\begin{center} \footnotesize $^{a}$
D$\acute{e}$pt. de math$\acute{e}$matiques, Universit$\acute{e}$ du
Qu$\acute{e}$bec $\grave{a}$ Montr$\acute{e}$al, C.P. 8888, Succ.
Centre-Ville Montr$\acute{e}$al, Canada H3C 3P8\\
$^{b,c}$ Institute of Mathematics, Academia Sinica, Taipei, Taiwan\\
 $^{d}$ Department of Mathematics, National Taiwan University,
Taipei, Taiwan
\end{center}
\date{}
 \vspace*{-0.3cm}
\thispagestyle{empty}
\begin{abstract}
The circular peak set of a permutation $\sigma$ is the set
$\{\sigma(i)\mid \sigma(i-1)<\sigma(i)>\sigma(i+1)\}$. In this
paper, we focus on the enumeration problems for permutations by
circular peak sets. Let $cp_n(S)$ denote the number of the
permutations of order $n$ which have the circular peak set $S$. For
the case with $|S|=0,1,2$, we derive the explicit formulas for
$cp_n(S)$. We also obtain some recurrence relations for the sequence
$cp_n(S)$ and give the formula for $cp_n(S)$ in the general case.
\end{abstract}
\noindent {\bf Keywords: Circular peak; Lattice path; Permutation; }

\newpage
\section{Introduction}
Throughout this paper, let $[m,n]:=\{m,m+1,\cdots,n\}$, $[n]:=[1,n]$
and $[m,n]=\emptyset$ if $m>n$. Let $\mathfrak{S}_n$ be the set of
all the permutations in the set $[n]$. We will write permutations of
$\mathfrak{S}_n$ in the form
$\sigma=(\sigma(1)\sigma(2)\cdots\sigma(n))$. We say that a
permutation $\sigma$ has a { circular descent} of value $\sigma(i)$
if $\sigma(i)>\sigma(i+1)$ for any $i\in[n-1]$. The { circular
descent set} of a permutation $\sigma$, denoted $CDES(\sigma)$, is
the set $\{\sigma(i)\mid \sigma(i)>\sigma(i+1)\}.$  For any
$S\subseteq [n]$, define a set $CDES_n(S)$ as
$CDES_n(S)=\{\sigma\in\mathfrak{S}_n\mid CDES(\sigma)=S\}$ and use
$cdes_n(S)$ to denote the number of the permutations in the set
$CDES_n(S)$, i.e., $cdes_n(S)=|CDES_n(S)|$. In a join work
\cite{chuang}, Hungyung Zhang et al. derive the explicit formula for
$cdes_n(S)$. As a application of the main results in \cite{chuang},
they also give the enumeration of permutation tableaux according to
their shape and generalizes the results in \cite{Dom}. Moreover,
Robert J.Clarke et al. \cite{clarke} gave the conceptions of linear
peak and cyclic peak and studied some new Mahonian permutation
statistics. In this paper, we say that a permutation $\sigma$ has a
{\it circular peak} of value $\sigma(i)$ if
$\sigma(i-1)<\sigma(i)>\sigma(i+1)$ for any $i\in[2,n-1]$. The {\it
circular peak set} of a permutation $\sigma$, denoted $CP(\sigma)$,
is the set $\{\sigma(i)\mid \sigma(i-1)<\sigma(i)>\sigma(i+1)\}.$
For example, the circular peak set of $\sigma=(48362517)$ is
$\{5,6,8\}$. Since $\sigma$ has no circular peaks when $n\leq 2$, we
always suppose that $n\geq 3$. For any $S\subseteq [n]$, we define a
set $CP_n(S)$ as $CP_n(S)=\{\sigma\in\mathfrak{S}_n\mid
CP(\sigma)=S\}.$ Obviously, if $\{1,2\}\subseteq S$, then
$CP_n(S)=\emptyset$.
\begin{exa} \begin{eqnarray*}CP_5(\{4,5\})=\{&14253,14352,24153,34152,24351,34251,&\\
&15243,15342,25143,35142,25341,35241&\}\end{eqnarray*}
\end{exa}
Suppose that $S=\{i_1,i_2,\cdots,i_k\}$, where $i_1<i_2<\cdots<i_k$.
Pierre Bouchard et al. \cite{Bouchard} proved that the necessary and
sufficient conditions for $CP_n(S)\neq \emptyset$ are $i_j\geq 2j+1$
for all $j\in [k]$.

Let $\mathcal{P}_n=\{S\mid CP_n(S)\neq\emptyset\}$. We can make the
set $\mathcal{P}_n$ into a poset $\mathscr{P}_n$ by defining
$S\preceq T$ if $S\subseteq T$ as sets. See \cite{Bouchard} for the
properties of the poset $\mathscr{P}_n$.

In this paper, we focus on the enumerations for permutations by
circular peak sets. Let $cp_n(S)$ denote the number of the
permutations in the sets $CP_n(S)$, i.e., $cp_n(S)=|CP_n(S)|$.  For
the case with $|S|=0,1,2$, we derive the explicit formulas for
$cp_n(S)$. For the general case, we consider the recurrence
relations for $cp_n(S)$. We find that if $n\geq 3$ and $S\subseteq
[n-1]$ with $CP_n(S)\neq \emptyset$, then $cp_n(S)$ satisfies the
following recurrence relation
$$cp_{n}(S\cup\{n\})=[n-2-2|S|]cp_{n-1}(S)+\sum\limits_{j\notin
S}2cp_{n-1}(S\cup\{j\}).$$ But we are more interested in another
recurrence relation as follows:\\
 suppose $n\geq 3$, $k\geq 0$ and
$S\subseteq [n-k-1]$ with $CP_n(S)\neq \emptyset$, then
\begin{eqnarray*}cp_{n}(S\cup[n-k+1,n])=2(k+1)cp_{n-1}(S\cup[n-k,n-1])+k(k+1)cp_{n-2}(S\cup[n-k,n-2]).\end{eqnarray*}
 First, let  $k\geq 1$, $n\geq \max\{3,2k\}$ and $CP_{n}([n-k+1,n])\neq\emptyset$, by this recurrence relation, we conclude that
$$cp_{n}([n-k+1,n])=k(k+1)\sum\limits_{i=1}^{k}
(-1)^{i+1}b_{k,i}\left[(2k+2)^{n-2k}-(2k+2-2i)^{n-2k}\right],$$
where the coefficients $b_{k,i}$ satisfy the following recurrence
relation:$$b_{k+1,i}=\left\{\begin{array}{lll}k(k+1)^2\sum\limits_{j=1}^k(-1)^{j+1}b_{k,j}&\text{\it
if}&i=1\\
k(k+1)\frac{k+2-i}{i}b_{k,i-1}&\text{\it if}&2\leq i\leq
k+1\end{array}\right.$$ with initial condition
$b_{1,1}=\frac{1}{2}$. And then we study the connections between
this recurrence relation and the circular-peak path with weight. A
{\it circular-peak path} is a lattice path in the first quadrant
starting at $(r,0)$ and ending at $(n,k)$ with only two kinds of
steps---\emph{horizon step $H=(1,0)$ }and \emph{rise step
$R=(2,1)$}. We consider a circular-peak path $P$ from $(r,0)$ to
$(n,k)$ as a word of $n-r-k$ letters using only $H$ and $R$. Let
$P_{r,n,k}$ be the set of all the circular-peak paths from the
vertices $(r,0)$ to $(n,k)$. Given an integer $i$ and
$P=e_1e_2\cdots e_{n-k-r}\in P_{r,n,k}$, if the step $e_j$ connects
the vertices $(x,y)$ with $(x+1,y)$, then the weight  of $e_j$,
denoted $w_i(e_j)$, is $2i+2(y+1)$; if the step $e_j$ connects the
vertices $(x,y)$ with $(x+2,y+1)$, then the weight
 of $e_j$, denoted $w_i(e_j)$, is $(y+i+1)(y+i+2)$; at last, let
$w_i(P)=\prod\limits_{j=1}^{n-k-r}w_i(e_j)$ be the weight of the
circular-peak path $P$ and  $w(i,r,n,k)=\sum\limits_{P\in
P_{r,n,k}}w_i(P)$. It is proved that if
  $k\geq 0$, $n\geq k+4$, $S\subseteq[3,n-k-1]$ with
$CP_n(S)\neq\emptyset$ and $r=\max S$, then
$$cp_{n}(S\cup[n-k+1,n])=\sum\limits_{i=0}^{k}w(i,r,n-i,k-i)cp_{r+i}(S\cup
[r+1,r+i]).$$ We must give the formula for $w(i,r,n,k)$ as follows:
$$w(i,r,n,k)=2^{n-r-2k}\prod\limits_{m=0}^{k-1}(m+i+1)(m+i+2)\sum\prod\limits_{m=0}^{k}[i+m+1]^{t_m},$$
where the sum is over all $(k+1)$-tuples $(t_0,t_1,\cdots,t_{k})$
such that $\sum\limits_{m=0}^{k}t_m=n-r-2k$ and $t_m\geq 0$. For any
$S\subseteq [3,n]$, define the {\it type} of the set $S$ as
$(r_1^{k_1},r_2^{k_2},\cdots,r_m^{k_m})$ if
$S=\bigcup\limits_{i=1}^{m}[r_i-k_i+1,r_i]$ such that $r_i\leq
r_{i+1}-k_{i+1}-1$ for all $i\in[1,m-1]$. We may state one of the
main results of the paper as follows: \begin{eqnarray*}
cp_{n}(S)&=&2^{n-r_m}\sum\limits_{i_1=0}^{k_m}\sum\limits_{i_2=0}^{k_{m-1}+i_1}\cdots\sum\limits_{i_{m-1}=0}^{k_{2}+i_{m-2}}\prod\limits_{j=1}^{m-1}\\
&&w(i_j,r_{m-j},r_{m-j+1}+i_{j-1}-i_{j},k_{m-j+1}+i_{j-1}-i_{j})cp_{r_1+i_{m-1}}([r_1-k_1+1,r_1+i_{m-1}]),
\end{eqnarray*}where $i_0=0$, when $m\geq 2$ and $CP_n(S)\neq\emptyset$.

The paper is organized as follows. In Section $2$, we will consider
the enumerations of the permutations in the sets $CP_n(S)$ with
$|S|=0,1,2$. In Section $3$, we will derive some recurrence
relations for the sequences $cp_n(S)$ and give the formula for
$cp_n(S)$ in the general case. In the Appendix, we list all the
values of $cp_n(S)>0$ for $3\leq n\leq 8$.
\section{The Enumerations for The Permutations in The Set $CP_n(S)$ with $|S|=0,1,2$}
In this section, we will consider the enumeration problems of the
permutations in the sets $CP_n(S)$ with $|S|=0,1,2$.

Let $cp_n(S)$ denote the number of the elements in the sets
$CP_n(\{S\})$, i.e., $cp_n(S)=|CP_n(S)|$. First, we need the
following lemma.
\begin{lem}\label{maxSn} Suppose $n\geq 3$ and
$S\subseteq [n]$ with $CP_n(S)\neq\emptyset$. Then

(1) $cp_{n+1}(S)=2cp_{n}(S)$, and

(2) let $m=\max S$, then $cp_n(S)=2^{n-m}cp_m(S)$ for any $n\geq m$.
\end{lem}
\begin{proof} (1) It is easy to check that
$((n+1)\sigma(1)\cdots\sigma(n))\in CP_{n+1}(S)$ and
$(\sigma(1)\cdots\sigma(n)(n+1))\in CP_{n+1}(S)$ for any
$\sigma=(\sigma(1)\cdots\sigma(n))\in CP_n(S)$. Conversely, for any
$\sigma\in CP_{n+1}(S)$, the position of the letter $n+1$ is $1$ or
$n+1$, i.e., $\sigma^{-1}(n+1)=1\text{ or }n+1$, since $n+1\notin
S$. Hence, $cp_{n+1}(S)=2cp_{n}(S)$.

(2) Iterating the identity in (1), we immediately obtain that
$cp_n(S)=2^{n-m}cp_m(S)$. \end{proof}

 For any
$\sigma\in\mathfrak{S}_n$, let $\tau$ be a subsequence
$(\sigma(j_1)\sigma(j_2)\cdots\sigma(j_k))$ of
$(\sigma(1)\cdots\sigma(n))$, where $1\leq j_1<j_2<\cdots <j_k\leq
n$. Define $red_{\sigma,\tau}$ as an increasing bijection of
$\{\sigma(j_1),\sigma(j_2),\cdots,\sigma(j_k)\}$ onto $[k]$. Let
$red_{\sigma}(\tau)=(red_{\sigma,\tau}(\sigma(j_1))red_{\sigma,\tau}(\sigma(j_2))\cdots
red_{\sigma,\tau}(\sigma(j_k)))$.
\begin{thm}\label{lemma|s|<=2} Let $n\geq 3$. Then
(1) $cp_n(\emptyset)=2^{n-1}$,

(2) $cp_n(\{i\})=2^{n-2}(2^{i-2}-1)$ for any $i\in[3,n]$, and

(3) $cp_n(\{i,j\})=2^{n-3}(2^{i-2}-1)(2^{j-i-1}-1)+2^{n+j-i-5}\cdot
3(3^{i-2}-2^{i-1}+1)$ for any $i,j\in[3,n]$ and $i<j$.
\end{thm}
\begin{proof} (1) For any $\sigma\in \mathfrak{S}_n$, suppose that
the position of the letter $1$ is $i+1$, i.e., $\sigma^{-1}(1)=i+1$.
Then $\sigma\in CP_n(\emptyset)$ if and only if $\sigma$ has the
form $\sigma(1)>\cdots
>\sigma(i+1)<\cdots <\sigma(n)$. For each letter $j\neq 1$, the position of $j$ has two possibilities at the left or
right of $1$. Hence, $cp_n(\emptyset)=2^{n-1}$.

(2) By Lemma \ref{maxSn}(2), we first consider the number of the
permutations in the set $CP_i(\{i\})$, where $i\geq 3$. For any
$\sigma\in CP_i(\{i\})$, suppose that the position of the letter $i$
is $k+1$, i.e., $\sigma^{-1}(i)=k+1$, then $1\leq k\leq i-2$,
$red_{\sigma}(\sigma(1)\cdots\sigma(k))\in CP_k(\emptyset)$ and
$red_{\sigma}(\sigma(k+2)\cdots\sigma(i))\in CP_{i-k-1}(\emptyset)$.
There are ${i-1\choose{k}}$ ways to form the set
$\{\sigma(1),\cdots,\sigma(k)\}$. So,
$cp_i(\{i\})=\sum\limits_{k=1}^{i-2}{i-1\choose{k}}2^{k-1}2^{i-k-2}=2^{i-2}(2^{i-2}-1)$.
Hence, $cp_n(\{i\})=2^{n-2}(2^{i-2}-1)$.

(3) It is easy to check that the identity holds when $i=3$ and
$j=4$. By Lemma \ref{maxSn}(2), we first consider the number of the
permutations in the set $CP_j(\{i,j\})$, where $3\leq i<j$. We begin
from the case $\sigma\in CP_j(\{i,j\})$ with
$\sigma^{-1}(i)<\sigma^{-1}(j)$. Let
\begin{eqnarray*}T_1(\sigma)&=&\{\sigma(k)\mid
\sigma(k)<i\text{ and }k<\sigma^{-1}(i)\},\\
T_2(\sigma)&=&\{\sigma(k)\mid
\sigma(k)<i\text{ and }\sigma^{-1}(i)<k<\sigma^{-1}(j)\},\\
T_3(\sigma)&=&\{\sigma(k)\mid \sigma(k)<i\text{ and
}k>\sigma^{-1}(j)\}.
\end{eqnarray*}
Note that $T_k(\sigma)\neq \emptyset$ for $k=1,2$ since $\sigma$
must have a circular peak $i$ and
$\bigcup\limits_{k=1}^3T_k(\sigma)=[i-1]$. Let
\begin{eqnarray*}T_4(\sigma)&=&\{\sigma(k)\mid
i<\sigma(k)<j\text{ and }k<\sigma^{-1}(i)\}\\
T_5(\sigma)&=&\{\sigma(k)\mid i<\sigma(k)<j\text{ and
}\sigma^{-1}(i)<k<\sigma^{-1}(j)\}.\end{eqnarray*} We discuss the
following two subcases.

{\it Subcase 1.} $T_3(\sigma)=\emptyset$

Let $T_6(\sigma)=\{\sigma(k)\mid i<\sigma(k)<j,k>\sigma^{-1}(j)\}$.
Then $T_6(\sigma)\neq\emptyset$ since $\sigma$ must have a circular
peak $j$ and $\bigcup\limits_{k=4}^6T_{k}(\sigma)=[i+1,j-1]$. For
$k=1,2,6$, the subsequences of $\sigma$, which is determined by the
elements in $T_{k}(\sigma)$, corresponds to a permutation in
$CP_{|T_{k}(\sigma)|}(\emptyset)$. The subsequences of $\sigma$,
which is determined by the elements in $T_{4}(\sigma)$ and
$T_{5}(\sigma)$, are decreasing and increasing, respectively. So,
the number of the permutations under this subcase is
\begin{eqnarray*}\sum\limits_{(T_1,T_2)}{i-1\choose{|T_1|,|T_2|}}2^{|T_1|-1}2^{|T_2|-1}
\sum\limits_{(T_4,T_5,T_6)}{j-i-1\choose{|T_4|,|T_5|,|T_6|}}2^{|T_6|-1}=2^{j-4}(2^{i-2}-1)(2^{j-i-1}-1),\end{eqnarray*}
where the first sum is over all pairs $(T_1,T_2)$ such that
$T_i\neq\emptyset$ for $i=1,2$ and $T_1\cup T_2=[i-1]$; the second
sum is over all triples $(T_4,T_5,T_6)$ such that $T_6\neq
\emptyset$ and $T_4\cup T_5\cup T_6=[i+1,j-1]$.

{\it Subcase 2.} $T_3(\sigma)\neq \emptyset$

Suppose that $\min T_3(\sigma)=s$, let
\begin{eqnarray*}T_6(\sigma)&=&\{\sigma(k)\mid
i<\sigma(k)<j\text{ and }\sigma^{-1}(j)<k<\sigma^{-1}(s)\}\\
T_7(\sigma)&=&\{\sigma(k)\mid i<\sigma(k)<j\text{ and
}k>\sigma^{-1}(s)\}.\end{eqnarray*}Then, for $k=1,2,3$, the
subsequences of $\sigma$, which is determined by the elements in
$T_{k}(\sigma)$, corresponds to a permutation in
$CP_{|T_{k}(\sigma)|}(\emptyset)$. The subsequences of $\sigma$,
which is determined by the elements in $T_{4}(\sigma)$ and
$T_{6}(\sigma)$, are decreasing. The subsequences of $\sigma$, which
is determined by the elements in $T_{5}(\sigma)$ and
$T_{7}(\sigma)$, are increasing. So, the number of the permutations
under this subcase is
\begin{eqnarray*}\sum\limits_{(T_1,T_2,T_3)}{i-1\choose{|T_1|,|T_2|,|T_3|}}2^{|T_1|-1}2^{|T_2|-1}2^{|T_3|-1}4^{j-i-1}=2^{2j-i-6}\cdot 3(3^{i-2}- 2^{i-1}+1)
,\end{eqnarray*} where the sum is over all triples $(T_1,T_2,T_3)$
such that $T_i\neq\emptyset$ for $i=1,2,3$ and $T_1\cup T_2\cup
T_3=[i-1]$.

Similarly, we may consider the case $\sigma\in CP_j(\{i,j\})$ with
$\sigma^{-1}(i)>\sigma^{-1}(j)$. Therefore,
$cp_{j}(\{i,j\})=2[2^{j-4}(2^{i-2}-1)(2^{j-i-1}-1)+2^{2j-i-6}\cdot
3(3^{i-2}- 2^{i-1}+1)]$. In general, for any $n\geq 3$ and $3\leq
i<j\leq n$,
$$cp_n(\{i,j\})=2^{n-3}(2^{i-2}-1)(2^{j-i-1}-1)+2^{n+j-i-5}\cdot
3(3^{i-2}-2^{i-1}+1).$$ \end{proof}

\section{The recurrence relations for the sequence $cp_n(S)$}
In this section, we will derive some recurrence relations for the
sequence $cp_n(S)$.

\begin{lem}\label{lemmarecurrencescupn}Let $n\geq 3$ and $S\subseteq[n-1]$. Then
\begin{eqnarray*}cp_{n}(S\cup\{n\})=[n-2-2|S|]cp_{n-1}(S)+\sum\limits_{j\notin S,j<n}2cp_{n-1}(S\cup\{j\}).\end{eqnarray*}
\end{lem}
\begin{proof} Suppose $\sigma\in CP_{n-1}(S)$.  We want to form a new
permutation $\tau\in CP_{n}(S\cup\{n\})$ by inserting the letter $n$
into $\sigma$. For any $j\in S$, since the letter $j$ still must be
a circular peak in the new permutation, we can't insert $n$  into
$\sigma$ beside $j$. But the letter $n$ must be a circular peak. So,
there are $(n-2-2|S|)$ ways to form a new permutation $\tau$ from
$\sigma$ such that $\tau\in CP_{n}(S\cup\{n\})$.

For any $j\notin S$ with $j<n$ and $\sigma\in CP_{n-1}(S\cup\{j\})$,
we must insert $n$  into $\sigma$ beside $j$ such that $n$ becomes a
circular peak. So, there are $2$ ways to form a new permutation
$\tau$ from $\sigma$ such that $\tau\in CP_{n}(S\cup\{n\})$.

Hence,
\begin{eqnarray*}cp_{n}(S\cup\{n\})=[n-2-2|S|]cp_{n-1}(S)+\sum\limits_{j\notin
S,j<n}2\cdot cp_{n-1}(S\cup\{j\}).\end{eqnarray*}
\end{proof}

For any $S\in [n]$, suppose $S=\{i_1,i_2,\ldots,i_k\}$, let
$\mathbf{x}_{S}$ stand for the monomial $x_{i_1}x_{i_2}\cdots
x_{i_k}$; In particular, let $\mathbf{x}_{\emptyset}=1$. Given
$n\geq 3$, we define a generating function
$${g_n(x_1,x_2,\ldots,x_{n};y) = \sum\limits_{\sigma \in
\mathfrak{S}_{n}} \mathbf{x}_{CP(\sigma)}}y^{|CP(\sigma)|}.$$ We
also write $g_n(x_1,x_2,\ldots,x_{n};y)$ as $g_n$ for short.

\begin{cor}Let $n$ be a positive integer with $n\geq 3$ and $g_n=\sum\limits_{\sigma \in
\mathfrak{S}_{n}} \mathbf{x}_{CP(\sigma)}y^{|CP(\sigma)|}.$ Then
$g_n$ satisfies the following recursion:
\begin{eqnarray*}g_{n+1}=[2+(n-1)x_{n+1}y]g_n+2x_{n+1}\sum\limits_{i=1}^{n}\frac{\partial g_n}{\partial
x_i}-2x_{n+1}y^2\frac{\partial g_n}{\partial y}.\end{eqnarray*} for
all $n\geq 3$ with initial condition $g_3=4+2x_3y$, where the
notation $\frac{\partial g_n}{\partial y}$ denote partial
differentiation of $g_n$ with respect to $y$.
\end{cor}
\begin{proof} Obviously, $g_3=4+2x_3y$ and  ${\sum\limits_{\sigma \in
\mathfrak{S}_{n}} \mathbf{x}_{CP(\sigma)}y^{|CP(\sigma)|} =
\sum\limits_{S\subseteq [2,n]}cp_{n}(S)\mathbf{x}_{S}}y^{|S|}.$
Hence,
\begin{eqnarray*}g_{n+1}
&=&\sum\limits_{S\subseteq [n+1]}cp_{n+1}(S)\mathbf{x}_{S}y^{|S|}\\
&=&\sum\limits_{S\subseteq [n+1],n+1\in
S}cp_{n+1}(S)\mathbf{x}_{S}y^{|S|}
+\sum\limits_{S\subseteq [n+1],n+1\notin S}cp_{n+1}(S)\mathbf{x}_{S}y^{|S|}\\
&=&\sum\limits_{S\subseteq
[n]}\left[(n-1-2|S|)cp_n(S)+\sum\limits_{i\in[n]\setminus
S}2cp_n(S\cup\{i\})\right]\mathbf{x}_{S}x_{n+1}y^{|S|+1}
+2g_n\\
&=&2\sum\limits_{S\subseteq [n]}\sum\limits_{i\in[n]\setminus
S}cp_n(S\cup\{i\})\mathbf{x}_{S}x_{n+1}y^{|S|+1}
-2\sum\limits_{S\subseteq
[n]}|S|cp_n(S)\mathbf{x}_{S}x_{n+1}y^{|S|+1}\\
&& +[2+(n-1)x_{n+1}y]g_n.
\end{eqnarray*}
Note that \begin{eqnarray*}\frac{\partial g_n}{\partial
y}=\sum\limits_{S\subseteq
[n]}|S|cp_n(S)\mathbf{x}_{S}y^{|S|-1}\end{eqnarray*} and
\begin{eqnarray*}&&\sum\limits_{S\subseteq [n]}\sum\limits_{i\in[n]\setminus
S}cp_n(S\cup\{i\})\mathbf{x}_{S}x_{n+1}y^{|S|+1}\\
&=&\sum\limits_{S\subseteq [n],S\neq
\emptyset}cp_n(S)x_{n+1}y^{|S|}\sum\limits_{i\in
S}\frac{\mathbf{x}_{S}}{x_{i}}\\
&=&x_{n+1}\sum\limits_{i=1}^{n}\frac{\partial g_n}{\partial
x_i}.\end{eqnarray*} Therefore,
\begin{eqnarray*}g_{n+1}=[2+(n-1)x_{n+1}y]g_n+2x_{n+1}\sum\limits_{i=1}^{n}\frac{\partial g_n}{\partial
x_i}-2x_{n+1}y^2\frac{\partial g_n}{\partial y}.\end{eqnarray*}
\end{proof}

\begin{lem}\label{lemmaScup[n-k+1,n]}Suppose that  $k\geq 0$ and $n\geq k+4$. Let
$S\subseteq[3,n-k-1]$ with $CP_n(S)\neq \emptyset$. Then
\begin{eqnarray*}cp_{n}(S\cup[n-k+1,n])=2(k+1)cp_{n-1}(S\cup[n-k,n-1])+k(k+1)cp_{n-2}(S\cup[n-k,n-2]).\end{eqnarray*}
\end{lem}
\begin{proof} For any $\sigma\in CP_{n}(S\cup[n-k+1,n])$, we consider
the following four cases.

{\it Case 1.} There is no letters $i\in [n-k+1,n]$ such that the
position of $i$ is beside $n-k$ in $\sigma$, i.e.,
$|\sigma^{-1}(i)-\sigma^{-1}(n-k)|=1$. Then
$\sigma^{-1}(n-k)=1\text{ or }n$ since the permutation $\sigma$
hasn't a circular peak $n-k$. We can obtain a new permutation $\tau$
by exchanging the positions of $n-k$  and $n$ in $\sigma$. Clearly,
$\tau\in CP_{n}(S\cup[n-k,n-1])$.  Lemma \ref{maxSn} (1) tells us
that $cp_{n}(S\cup[n-k,n-1])=2cp_{n-1}(S\cup[n-k,n-1])$. Hence, the
number of permutations under this case is $2\cdot
cp_{n-1}(S\cup[n-k,n-1])$.

{\it Case 2.} There are exact two letters $j,m\in [n-k+1,n]$ such
that $|\sigma^{-1}(j)-\sigma^{-1}(n-k)|=1$ and
$|\sigma^{-1}(m)-\sigma^{-1}(n-k)|=1$.  Deleting $j$ and $m$, we
obtain a subsequence $\tau$ of $\sigma$ . Then
$red_{\sigma}(\tau)\in CP_{n-2}(S\cup[n-k,n-2])$. Note that there
are $k(k-1)$ ways to form the pairs $(j,m)$. Hence, the number of
permutations under this case is $k(k-1)cp_{n-2}(S\cup[n-k,n-2])$.

{\it Case 3.} There is exact one letter $j\in [n-k+1,n]$ such that
$|\sigma^{-1}(j)-\sigma^{-1}(n-k)|=1$. Then there are $k$ ways to
form the set $\{j\}$. Let $\tau$ be the subsequence of $\sigma$
obtained by deleting $j$. There are the following two subcases.

 {\it Subcase 3.1.}
$\sigma^{-1}(n-k)\neq 1\text{ and }n$.  Then $red_{\sigma}(\tau)\in
CP_{n-1}(S\cup[n-k,n-1])$. Hence, the number of permutations under
this subcase is $k\cdot cp_{n-1}(S\cup[n-k,n-1])$.

{\it Subcase 3.2.} $\sigma^{-1}(n-k)= 1\text{ or }n$. Then
$red_{\sigma}(\tau)\in CP_{n-2}(S\cup[n-k,n-2])$. Hence, the number
of permutations under this subcase is $k\cdot
cp_{n-2}(S\cup[n-k,n-2])$.

So,
\begin{eqnarray*}&&cp_{n}(S\cup[n-k+1,n])\\
&=&2cp_{n-1}(S\cup[n-k,n-1])+k(k-1)cp
_{n-2}(S\cup[n-k,n-2])\\
&&+2k\cdot cp_{n-1}(S\cup[n-k,n-1])+2k\cdot cp_{n-2}(S\cup[n-k,n-2])\\
&=&2(k+1)cp_{n-1}(S\cup[n-k,n-1])+k(k+1)cp_{n-2}(S\cup[n-k,n-2]).\end{eqnarray*}
\end{proof}

Setting $S=\emptyset$ in Lemma \ref{lemmaScup[n-k+1,n]}, we derive
the following results.

\begin{thm}\label{thmcoeb} Let  $k\geq 1$, $n\geq max\{3,2k\}$ and $CP_{n}([n-k+1,n])\neq\emptyset$. Suppose that
$$cp_{n}([n-k+1,n])=k(k+1)\sum\limits_{i=1}^{k}
(-1)^{i+1}b_{k,i}\left[(2k+2)^{n-2k}-(2k+2-2i)^{n-2k}\right].$$ Then
the coefficients $b_{k,i}$ satisfy the recurrence relation as
follows:$$b_{k+1,i}=\left\{\begin{array}{lll}k(k+1)^2\sum\limits_{j=1}^k(-1)^{j+1}b_{k,j}&\text{\it
if}&i=1\\
k(k+1)\frac{k+2-i}{i}b_{k,i-1}&\text{\it if}&2\leq i\leq
k+1\end{array}\right..$$
\end{thm}
\begin{proof} Lemma \ref{lemma|s|<=2} tells us that
$cp_{n}([n,n])=cp_{n}(\{n\})=4^{n-2}-2^{n-2}$. Hence,
$b_{1,1}=\frac{1}{2}$. Lemma \ref{lemmaScup[n-k+1,n]} implies that
\begin{eqnarray*}cp_{n}([n-k+1,n])=2(k+1)cp_{n-1}([n-k,n-1])+k(k+1)cp_{n-2}([n-k,n-2]).\end{eqnarray*}
By comparing the coefficients, we obtain the desired
results.\end{proof}

By Theorem \ref{thmcoeb}, we compute the values of $b_{k,i}$ for all
$1\leq k\leq 4$ as follows. {
$$
\begin{array}{|r|l|l|l|l|}
\hline
&i=1&2&3&4\\
 \hline
k=1&\frac{1}{2}&&&\\
 \hline
2&2&\frac{1}{2}&&\\
 \hline
3&27&12&1&\\
 \hline
 4&768&486&96&3\\
 \hline
\end{array}
$$}\begin{center} Table.1. The values of $b_{k,i}$ for $1\leq k\leq 4$\end{center}

\begin{cor}\label{corcoea} Let  $k\geq 0$, $n\geq max\{3,2k\}$ and $CP_{n}([n-k+1,n])\neq\emptyset$. Suppose that $$cp_{n}([n-k+1,n])=\sum\limits_{i=0}^{k}
(-1)^{i}a_{k,i}(2k+2-2i)^{n-2k}.$$ Then the coefficients $a_{k,i}$
satisfy the following recurrence
relation:$$a_{k+1,i}=\left\{\begin{array}{lll}\sum\limits_{j=1}^{k+1}(-1)^{j+1}a_{k+1,j}&\text{\it
if}&i=0\\
(k+1)(k+2)\frac{k+2-i}{i}a_{k,i-1}&\text{\it if}&1\leq i\leq
k+1\end{array}\right.,$$ with initial condition
$a_{0,0}=\frac{1}{2}$. Let $f_k(x)=\sum\limits_{i=0}^ka_{k,i}x^i$
for any $k\geq 0$, then $f_0(x)=\frac{1}{2}$, $f_{k+1}(-1)=0$ and
$f_k(x)$ satisfies the differential equation
$f_{k+1}'(x)=(k+1)^2(k+2)f_{k}(x)-(k+1)(k+2)xf_{k}'(x)$ for any
$k\geq 0$, where the notation ``$~'~$" denotes the differentiation
of functions.
\end{cor}
\begin{proof}  When $k=0$, by Lemma \ref{lemma|s|<=2}, we have
$cp_{n}([n+1,n])=cp_{n}(\emptyset)=2^{n-1}$. Hence,
$a_{0,0}=\frac{1}{2}$. When $k\geq 1$, it is easy to check that
$$a_{k,i}=\left\{\begin{array}{lll}k(k+1)\sum\limits_{j=1}^k(-1)^{j+1}b_{k,j}&\text{\it
if}&i=0\\
k(k+1)b_{k,i}&\text{\it if}&1\leq i\leq k\end{array}\right..$$
Hence, by Theorem \ref{thmcoeb}, we have
$$a_{k+1,i}=\left\{\begin{array}{lll}\sum\limits_{j=1}^{k+1}(-1)^{j+1}a_{k+1,j}&\text{\it
if}&i=0\\
(k+1)(k+2)\frac{k+2-i}{i}a_{k,i-1}&\text{\it if}&1\leq i\leq
k+1\end{array}\right.,$$ with initial condition
$a_{0,0}=\frac{1}{2}$. Let $f_k(x)=\sum\limits_{i=0}^ka_{k,i}x^i$
for any $k\geq 0$. Clearly, $f_0(x)=\frac{1}{2}$. Note that
$a_{k,i}=\sum\limits_{j=1}^{k}(-1)^{j+1}a_{k,j}$ for any $k\geq 1$.
Hence, $f_{k+1}(-1)=0$ if $k\geq 0$. Since
$a_{k+1,i}=(k+1)(k+2)\frac{k+2-i}{i}a_{k,i-1}$ for any $i\in [k+1]$,
we have
\begin{eqnarray*}\sum\limits_{i=1}^{k+1}a_{k+1,i}x^i=\sum\limits_{i=1}^{k+1}(k+1)(k+2)\frac{k+2-i}{i}a_{k,i-1}x^i\end{eqnarray*}
Simple computations tell us that
$f_{k+1}'(x)=(k+1)^2(k+2)f_{k}(x)-(k+1)(k+2)xf_{k}'(x)$.\end{proof}

By Corollary \ref{corcoea}, we compute the values of $a_{k,i}$ for
all $0\leq k\leq 3$ as follows. {
$$
\begin{array}{|r|l|l|l|l|}
\hline
&i=0&1&2&3\\
 \hline
k=0&\frac{1}{2}&&&\\
 \hline
1&1&1&&\\
 \hline
2&9&12&3&\\
 \hline
 3&192&324&144&12\\
 \hline
\end{array}
$$}\begin{center} Table.2. The values of $a_{k,i}$ for $0\leq i\leq 3$\end{center}

Recall that $w(i,r,n,k)=\sum\limits_{P\in P_{r,n,k}}w_i(P)$, where
$P_{r,n,k}$ is the set of all the circular-peak path from the
vertices $(r,0)$ to $(n,k)$ and $w_i(P)$ is the weight of path $P\in
P_{r,n,k}$.

\begin{lem}\label{lemmaweigth}\begin{eqnarray*}w(i,r,n,k)=2^{n-r-2k}\prod\limits_{m=0}^{k-1}(m+i+1)(m+i+2)\sum\prod\limits_{m=0}^{k}(i+m+1)^{t_m}\end{eqnarray*}
where the sum is over all $(k+1)$-tuples $(t_0,t_1,\cdots,t_{k})$
such that $\sum\limits_{m=0}^{k}t_m=n-r-2k$ and $t_m\geq 0$.
\end{lem}
\begin{proof} Suppose that $P=e_1e_2\cdots e_{n-k-r}\in P_{r,n,k}$.
Let $\mathcal {R}=\{j\mid e_j=R\}$, then $|\mathcal {R}|=k$.
Furthermore, we may suppose that $\mathcal {R}=\{e_{j_1},\cdots,
e_{j_k}\}$, where $0=j_0< j_1<j_2<\cdots <j_k\leq n-k-r=j_{k+1}$.
Hence,
$$w_i(P)=\prod\limits_{m=0}^{k}[2i+2m+2]^{j_{m+1}-j_m-1}\prod\limits_{m=0}^{k-1}(m+i+1)(m+i+2).$$
Let $t_m=j_{m+1}-j_m-1$ for any $0\leq m\leq k$, then $t_m\geq 0$
and $\sum\limits_{m=0}^{k}t_m=n-r-2k$. So,
\begin{eqnarray*}w(i,r,n,k)&=&\sum\prod\limits_{m=0}^{k}[2i+2m+2]^{t_m}\prod\limits_{m=0}^{k-1}(m+i+1)(m+i+2)\\
&=&2^{n-r-2k}\prod\limits_{m=0}^{k-1}(m+i+1)(m+i+2)\sum\prod\limits_{m=0}^{k}[i+m+1]^{t_m}\end{eqnarray*}
where the sum is over all $(k+1)$-tuples $(t_0,t_1,\cdots,t_{k})$
such that $\sum\limits_{m=0}^{k}t_m=n-r-2k$ and $t_m\geq
0$.\end{proof}
\begin{lem}\label{coefficient}Suppose $n\geq 3$ and $k\geq 0$. Let
$S\subseteq[3,n-k-1]$ with $CP_n(S)\neq\emptyset$ and $r=max S$.
Then
\begin{eqnarray*}cp_{n}(S\cup[n-k+1,n])=\sum\limits_{i=0}^{k}w(i,r,n-i,k-i)cp_{r+i}(S\cup [r+1,r+i]).\end{eqnarray*}
\end{lem}
\begin{proof} We view the set $S\cup [n-k+1,n]$ as a vertex $(n,k)$.
Connect the vertices $(n,k)$ with $(n-1,k)$ $( resp. (n-2,k-1))$ and
give this edge a weight $2(k+1)$ $( resp. k(k+1))$. We draw the
obtained graph as follows:
\begin{center}
\includegraphics[width=2in]{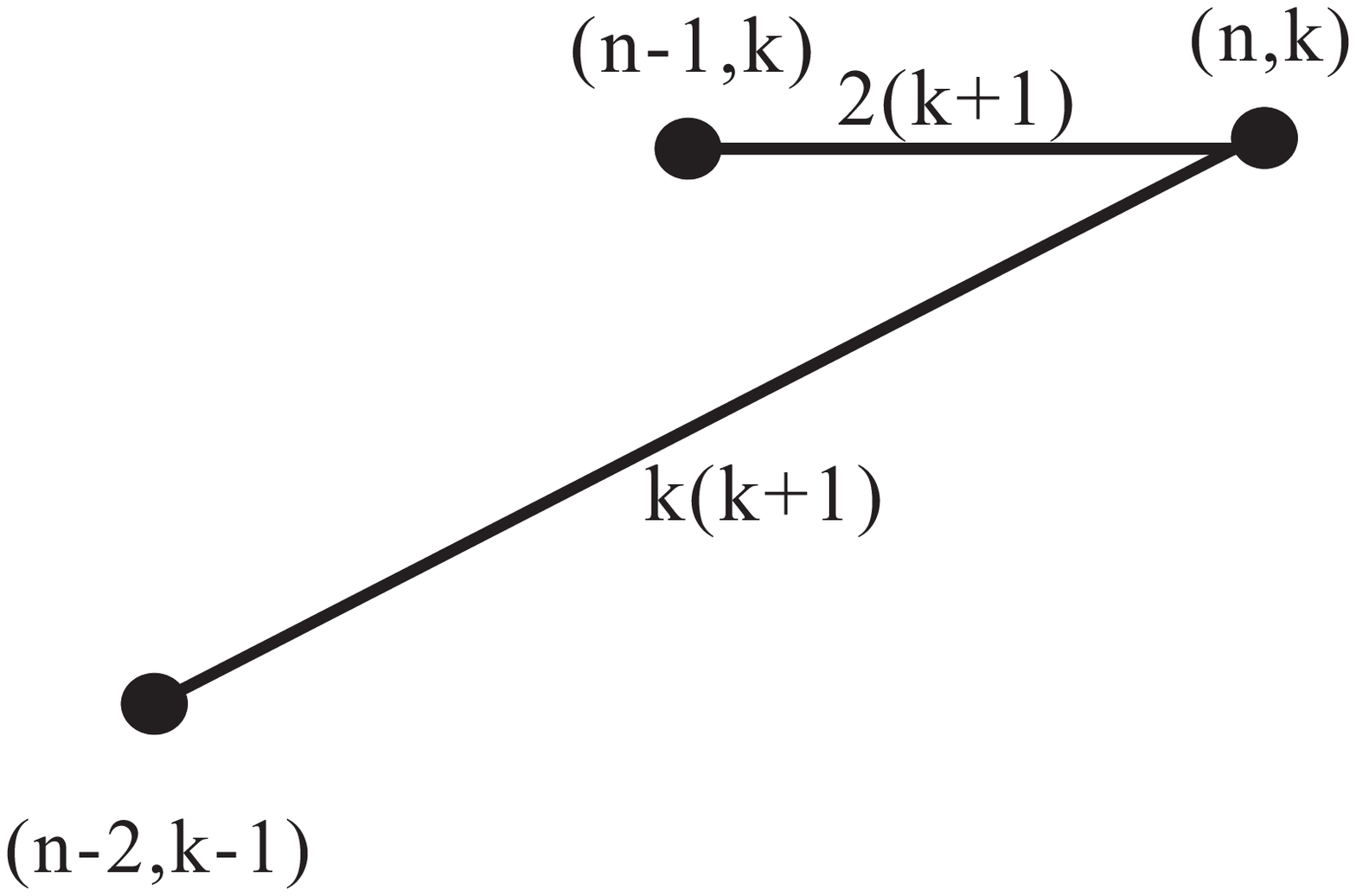}\\
 Fig.1. the obtained graph from the recurrence relation in Lemma
 \ref{lemmaScup[n-k+1,n]}
\end{center}
Conversely, from the graph, we also can derive the recurrence
relation in Lemma \ref{lemmaScup[n-k+1,n]}. Repeating this
processes,  we will obtain the following graph with weight.
\begin{center}
\includegraphics[width=4in]{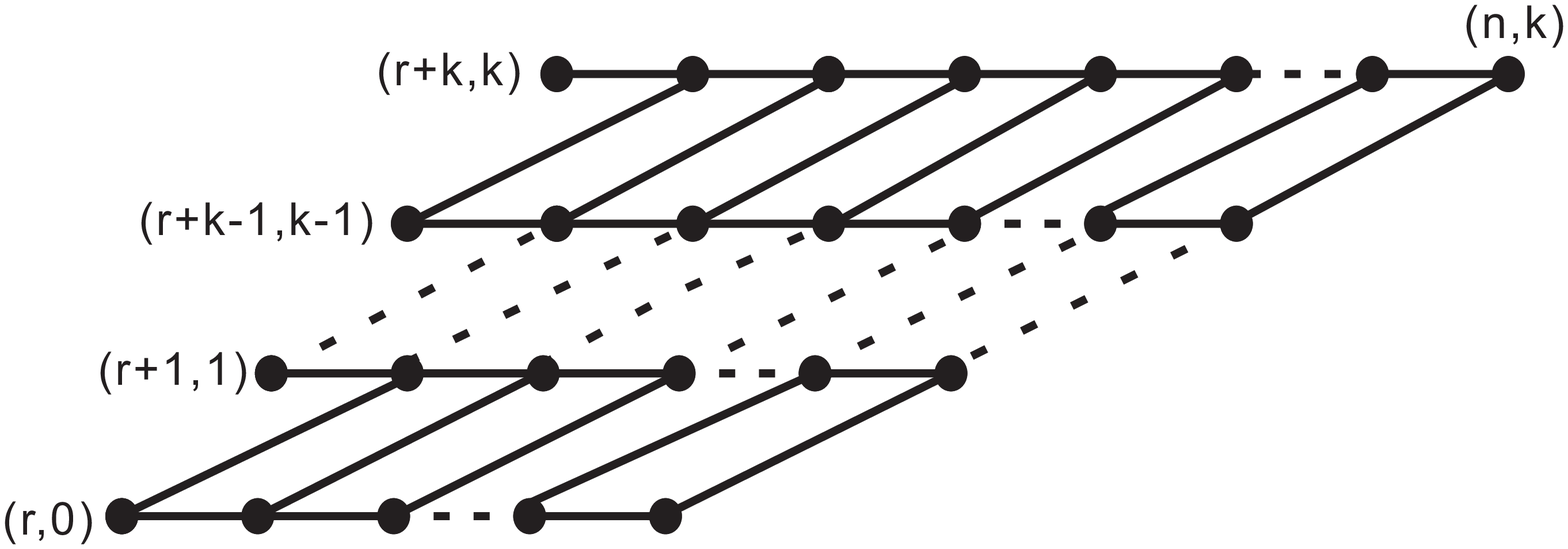}
\end{center}
 Hence, $cp_{n}(S\cup[n-k+1,n])$ can  be expressed as a
linear combination of $cp_{r+i}(S\cup [r+1,r])$, $cp_{r+i}(S\cup
[r+1,r+1])$, $\cdots$, and $cp_{r+i}(S\cup [r+1,r+k])$.  It is easy
to check that  the coefficient of $cp_{r+i}(S\cup [r+1,r+i])$ is
$w(i,r,n-i,k-i)$ which is the sum of the weights of all the circular
peak path from $(r+i,i)$ to $(n,k)$. So,
$cp_{n}(S\cup[n-k+1,n])=\sum\limits_{i=0}^{k}w(i,r,n-i,k-i)cp_{r+i}(S\cup
[r+1,r+i]).$\end{proof}

For any $S\subseteq [3,n]$, suppose that the type of the set $S$ is
$(r_1^{k_1},r_2^{k_2},\cdots,r_m^{k_m})$. We have studied the case
with $m=1$ in Theorem \ref{thmcoeb}. We will consider the case with
$m\geq 2$ in the following theorem.

\begin{thm}\label{theoremlastformular}Let $n\geq 3$, $k\geq 0$ and
$S\subseteq[3,n]$. Suppose that $CP_n(S)\neq \emptyset$ and the type
of the set $S$ is $(r_1^{k_1},r_2^{k_2},\cdots,r_m^{k_m})$ with
$m\geq 2$, where $m\geq 2$. Then
\begin{eqnarray*}
cp_{n}(S)&=&2^{n-r_m}\sum\limits_{i_1=0}^{k_m}\sum\limits_{i_2=0}^{k_{m-1}+i_1}\cdots\sum\limits_{i_{m-1}=0}^{k_{2}+i_{m-2}}\prod\limits_{j=1}^{m-1}\\
&&w(i_j,r_{m-j},r_{m-j+1}+i_{j-1}-i_{j},k_{m-j+1}+i_{j-1}-i_{j})cp_{r_1+i_{m-1}}([r_1-k_1+1,r_1+i_{m-1}])
\end{eqnarray*}where $i_0=0$.
\end{thm}
\begin{proof} Iterating the identity in Lemma \ref{coefficient}, we
obtain the desired results.\end{proof}

\begin{exa} Let $S=\{i,j\}$ such that $3\leq i<j-1$, $j\leq n$ and
$CP_n(S)\neq\emptyset$. The type of $S$ is $(i^1,j^1)$. By Theorem
\ref{theoremlastformular}, we have
$cp_n(\{i,j\})=2^{n-j}[w(0,i,j,1)cp_i(\{i\})+w(1,i,j-1,0)cp_{i+1}(\{i,i+1\})]$.
Lemma \ref{lemmaweigth} implies that
$w(0,i,j,1)=2^{j-i-1}(2^{j-i-1}-1)$ and $w(1,i,j-1,0)=2^{2j-2i-2}$.
Theorem \ref{lemma|s|<=2} tells us that
$cp_{i}(\{i\})=2^{i-1}(2^{i-1}-1)$. From Corollary \ref{corcoea} it
follows that
$cp_{i+1}(\{i,i+1\})=\sum\limits_{m=0}^{2}(-1)^{m}a_{2,m}(6-2m)^{i-3}=9\cdot
6^{i-3}-12\cdot 4^{i-3}+3\cdot 2^{i-3}.$ Hence, we obtain the
formula
$cp_n(\{i,j\})=2^{n-3}(2^{i-2}-1)(2^{j-i-1}-1)+2^{n+j-i-5}\cdot
3(3^{i-2}-2^{i-1}+1)$ again.
\end{exa}

\section{Appendix}
For convenience to check the identities given in the previous
sections, by the computer search, for $3\leq n\leq 8$, we obtain the
number $cp_n(S)$ of the permutations in the sets
$CP_n(S)\neq\emptyset$ and list them in Table $3$.

{\tiny $$
\begin{array}{|r|l|l|l|l|l|l|l|l|l|}
\hline
 n=3&S=\emptyset&\{3\}&&&&&&&\\
 \hline
&4&2&&&&&&&\\
 \hline
n=4&S=\emptyset&\{3\}&\{4\}&&&&&&\\
 \hline
&8&4&12&&&&&&\\
 \hline
n=5&S=\emptyset&\{3\}&\{4\}&\{5\}&\{3,5\}&\{4,5\}&&&\\
 \hline
&16&8&24&56&4&12&&&\\
 \hline
n=6&S=\emptyset&\{3\}&\{4\}&\{5\}&\{6\}&\{3,5\}&\{3,6\}&\{4,5\}&\{4,6\}\\
 \hline
&32&16&48&112&240&8&24&24&72\\
 \hline
&S=\{5,6\}&&&&&&&&\\
 \hline
 &144&&&&&&&&\\
 \hline
 n=7&S=\emptyset&\{3\}&\{4\}&\{5\}&\{6\}&\{7\}&\{3,5\}&\{3,6\}&\{3,7\}\\
 \hline
&64&32&96&224&480&992&16&48&112\\
 \hline
 &S=\{4,5\}&\{4,6\}&\{4,7\}&\{5,6\}&\{5,7\}&\{6,7\}&\{3,5,7\}&\{3,6,7\}&\{4,5,7\}\\
 \hline
&48&144&336&288&688&1200&8&24&24\\
 \hline
 &S=\{4,6,7\}&\{5,6,7\}&&&&&&&\\
 \hline
 &72&144&&&&&&&\\
 \hline
 n=8&S=\emptyset&\{3\}&\{4\}&\{5\}&\{6\}&\{7\}&\{8\}&\{3,5\}&\{3,6\}\\
 \hline
 &128&64&192&448&960&1984&4032&32&96\\
 \hline
 &S=\{3,7\}&\{3,8\}&\{4,5\}&\{4,6\}&\{4,7\}&\{4,8\}&\{5,6\}&\{5,7\}&\{5,8\}\\
 \hline
 &224&480&96&288&672&1440&576&1376&2976\\
 \hline
 &S=\{6,7\}&\{6,8\}&\{7,8\}&\{3,5,7\}&\{3,5,8\}&\{3,6,7\}&\{3,6,8\}&\{3,7,8\}&\{4,5,7\}\\
 \hline
&2400&5280&8640&16&48&48&144&288&48\\
 \hline
 &S=\{4,5,8\}&\{4,6,7\}&\{4,6,8\}&\{4,7,8\}&\{5,6,7\}&\{5,6,8\}&\{5,7,8\}&\{6,7,8\}&\\
 \hline
&144&144&432&864&288&864&1728&2880&\\
 \hline
\end{array}
$$}\begin{center} Table.3. $cp_{n}(S)$ for $3\leq n\leq 8$ with $CP_n(S)\neq\emptyset$\end{center}



\end{document}